\documentclass[12pt,a4paper]{article}
\usepackage[cp1251]{inputenc}
\usepackage[english,russian]{babel}
\usepackage{graphicx}
\usepackage{amssymb}
\usepackage{amsmath}
\usepackage{amsthm}
\usepackage{euscript}

\newtheorem{definition}{Определение}
\newtheorem{theorem}{Теорема}
\newtheorem{lemma}{Лемма}
\newtheorem{rem}{Замечание}
\newenvironment{remark}{\begin{rem}\rm}{\end{rem}}

\newcommand{\RP}[0]{\mathbb{RP}}
\newcommand{\CP}[0]{\mathbb{CP}}
\newcommand{\R}[0]{\mathbb{R}}
\newcommand{\Z}[0]{\mathbb{Z}}
\newcommand{\A}[0]{\mathfrak A}

\newcommand{\be}[0]{\begin{equation}}
\newcommand{\ee}[0]{\end{equation}}
\newcommand{\bez}[0]{\begin{equation*}}
\newcommand{\eez}[0]{\end{equation*}}
\newcommand{\bl}[0]{\begin{lemma}}
\newcommand{\el}[0]{\end{lemma}}
\newcommand{\lra}[0]{\Leftrightarrow}
\newcommand{\ep}[0]{\begin{flushright} $\square$ \end{flushright}}

\newcommand{\paragraf}[1]{
\medskip{\centerline{\bf #1}}\medskip}
\newcommand{\prove}[1]{
\smallskip \noindent{\it Доказательство #1.}~}

\begin{document}

{\centerline {\bf On the number of regions and multiplicities of vertices in plane arrangements}}

\medskip
{\small
For an arrangement of $n$ pseudolines in the real projective plane let us denote by $t_i$ the number of vertices incident to $i$ lines. Then
\bez
t_2+\frac 32 t_3\geq 8+\sum_{i\geq 4}\left(2i-7\frac 12\right)t_i
\eez
provided $t_n=t_{n-1}=t_{n-2}=0$. This looks similar the Hirzebruch inequality for $t_i$ numbers but has an elementary proof. We present an algorithm to get lower bounds of the number of regions basing on linear on $t_i$ inequalities like the above-mentioned. The lower bounds arised in connection with Martinov theorem on the set of all possible numbers of regions and we show how the new bounds may be applied in it.

}

\section{Введение}

\bigskip
На вещественной проективной плоскости $\RP^2$ будем рассматривать конечные наборы\footnote{термин {\it arrangement of lines} здесь переводится как {\it набор прямых}. Есть схожий объект: {\it configurations}~---~{\it конфигурации},  наборы точек и прямых с дополнительными требованиями на число инциденций.}, состоящие из $n\geq 2$ различных прямых.  Набору прямых соответствует представление плоскости $\RP^2$ в виде клеточного комплекса $\A$ с $v$ вершинами, $e$ ребрами (отрезками) и $f$ двумерными клетками. Вершины комплекса --- точки пересечения прямых, ребра --- отрезки прямых без внутренних точек пересечения, двумерные клетки --- замыкания компонент связности дополнения к объединению прямых.
Инвариантом наборов прямых является массив чисел $t_i$, где $t_i$ равно числу точек пересечения, принадлежащих $i$ прямым, $2 \leq i \leq n.$ Из известных результатов о наборах прямых на плоскости $\RP^2$ или о сходных объектах (см. обзоры \cite{Brass}, \cite{Nilakantan 2005}) нас будут интересовать следующие.

\begin{itemize}
  \item Нижние оценки числа двумерных клеток и описание всех возможных значений этого числа (\cite{Cordovil 80}, \cite{Grunbaum 72}, \cite{Martinov 90}, \cite{Martinov 93}, \cite{Purdy 80}, \cite{Arnold 08} и др.)
  \item Соотношения между числами прямых, вершин, ребер и двумерных клеток комплекса $\A$, между числами $t_i$ (\cite{Csima Sawyer 93}, \cite{Csima Sawyer 95}, \cite{Erdos Purdy 78}, \cite{Grunbaum 72}, \cite{Hirzebruch  86}, \cite{Melchior 40} и др.).
\end{itemize}

 После обзора результатов мы подробно обсудим теорему Мартинова \cite{Martinov 93} о классификации всех возможных чисел областей и ее следствия.
Для обоснования теоремы Мартинова понадобились нижние оценки числа областей при фиксированном максимальном числе прямых, проходящих через одну точку. Серию таких оценок мы докажем в разделе \ref{subsection ocenki}. В свою очередь, эффективные оценки можно строить с помощью линейных по $t_i$ неравенств. В разделе \ref{subsection linear t_i ine} мы перечислим известные и получим еще одно линейное по $t_i$ неравенство типа неравенства Хирцебруха \cite{Hirzebruch  86}.

Если у всех прямых набора существует общая точка, то набор прямых называется {\it тривиальным}. В дальнейшем наборы прямых предполагаются нетривиальными (если не сказано обратное).
{\it Изоморфными} называются наборы, клеточные комплексы которых изоморфны, т.е. для которых существует взаимно-однозначное соответствие между вершинами, ребрами и 2-клетками, сохраняющее все инциденции. Компоненты связности дополнения в плоскости $\RP^2$ к объединению прямых будем называть {\it областями}. Замыкания областей для нетривиальных наборов суть многоугольники, число $j$--угольников среди которых обозначим через $p_j$. Числа $p_j$ наряду с числами $t_i$ являются инвариантами наборов прямых при изоморфизме.

Многие из приведенных ниже фактов верны не только для наборов прямых на плоскости $RP^2$, но и для наборов кривых с некоторыми топологическими свойствами, присущими наборам прямых.
Назовем {\it псевдопрямой}\footnote{ перевод {\it pseudoline}} гладкую замкнутую несамопересекающуюся  кривую на вещественной проективной плоскости, не гомотопную отображению окружности в точку.

\begin{definition} Под набором псевдопрямых будем иметь в виду конечный набор из $n \geq 3$ псевдопрямых, любые две из которых пересекаются в единственной точке и пересекаются трансверсально.
\end{definition}
Нетривиальным набором псевдопрямых будем называть набор, в котором не все псевдопрямые проходят через одну точку. Будем рассматривать только нетривиальные наборы псевдопрямых.
Заметим, что не для всякого набора псевдопрямых существует изоморфный набор прямых (см. спрямляемость наборов псевдопрямых в обзоре \cite{Nilakantan 2005}). Пример такого набора строится с помощью теоремы Дезарга о том, что точки пересечения соответствующих сторон двух перспективных треугольников лежат на одной прямой. Треугольники называются {\it перспективными}, если три прямые, проходящие через соответствующие стороны треугольников, пересекаются в одной точке. Таким образом, неспрямляемым набором псевдопрямых будет набор из 6 сторон треугольников, трех прямых, проходящих через соответствующие вершины треугольников и псевдопрямой, проходящей через две точки пересечения соответствующих сторон треугольников и не проходящая через третью точку пересечения.

Имеется связь наборов прямых с зоноэдрами \cite{Grunbaum 67}. Зоноэдром в $\R^d$ называется цент\-раль\-но-сим\-мет\-рич\-ный выпуклый многогранник, все гиперграни которого центрально симметричны. Набор прямых в $\RP^2$ представим как набор проходящих через начало координат плоскостей в $\R^3$. Для каждой плоскости возьмем отрезок, состоящий из ее нормалей в начале координат. Сумма Минковского этих отрезков будет зоноэдром в $\R^3$, соответствующим набору прямых. При этом прямые соответствуют пояскам зоноэдра, области плоскости $\RP^2$ --- парам противоположных вершин зоноэдра, точки пересечения прямых --- парам противоположных граней зоноэдра. Например, трем не коллинеарным прямым на плоскости $\RP^2$ соответствует куб.

\paragraf{Соотношения между числами $n, v, e, f, t_i$ и $p_j$.}

Из формулы для эйлеровой характеристики следует
\bez
v-e+f=1
\eez
Подсчитывая число пар прямых, получаем
\begin{equation*}
\sum_{i\geq 2}i(i-1)t_i=n(n-1).
\end{equation*}
Число областей $f$ оценивается по индукции по $n$:
\begin{equation*}
2n-2 \leq f \leq 1+\frac{n(n-1)}2,
\end{equation*}
причем правое неравенство обращается в равенство для наборов псевдопрямых общего положения, а левое для наборов, в которых $n-1$ псевдопрямых проходят через одну точку.

Псевдопрямые на плоскости $\RP^2$ находятся в {\it общем положении}, если никакие три из них не имеют общей точки. { \it Симплициальными} называются нетривиальные наборы, для которых каждая область (компонента связности дополнения) примыкает по дугам ровно к трем псевдопрямым. Обозначим через $m$ максимальное число псевдопрямых набора, имеющих общую точку (т.е. $m = \max\{i \ | \ t_i \neq 0\}$).

\begin{theorem} {\rm Б.~Грюнбаум, \cite[pp. 401-402]{Grunbaum 67}} Для нетривиальных наборов прямых
\begin{equation*}
v+1 \leq f \leq 2v-2
\end{equation*}
причем неравенство слева обращается в равенство тогда и только тогда, когда прямые в общем положении, а неравенство справа обращается в равенство тогда и только тогда, когда набор прямых симплициальный.
\end{theorem}

Неравенство Мельхиора \cite{Melchior 40}
\bez
t_2\geq 3+\sum_{i\geq 4}(i-3)t_i
\eez
обращается в равенство для симплициальных наборов и, по-видимому, Мельхиор первым предпринял частичное описание симплициальных наборов псевдопрямых.

Из неравенства Мельхиора следует $\max\{t_2,t_3\} \geq t_i$ для всех $i\geq 2$. П.~Эрдош и Г.Б.~Пурди в \cite{Erdos Purdy 78} доказали, что
\begin{equation*}
 \max\{t_2,t_3\}\geq n-1 \quad \text{при} \quad n\geq 25.
\end{equation*}
Также они доказали, что если $t_2<n-1$, то
 $t_3>cn^2$ для некоторой положительной константы $c$.

Дж.~Сцима и Е.Т.~Сойер в \cite{Csima Sawyer 93, Csima Sawyer 95}  доказали, что $t_2\geq \frac{6}{13}n$ при $n\geq 8$. Известен пример с $t_2=\frac n2$ для четных $n\geq 6$, и гипотеза Дирака утверждает, что $t_2\geq \left[\frac n2\right]$.

Числа вершин, ребер и областей  можно выразить через $p_j$ следующим образом
\bez
v=1+ \frac 12 \sum_{j\geq 2} (j-2)p_j,\qquad
e= \frac 12 \sum_{j\geq 2} jp_j, \qquad
 f = \sum_{j\geq 2} p_j.
\eez

\paragraf{Оценки и возможные значения числа областей}

Грюнбаум в \cite{Grunbaum 72} впервые поставил вопрос об описании множества возможных чисел областей и доказал, что
\begin{equation*}
f \geq 3n-6 \quad \text{при} \quad m\leq n-2.
\end{equation*}
Тем самым, число областей не может принадлежать интервалу $(2n-2; 3n-6)$.
Грюнбаум предположил
\cite[conj. 2.4]{Grunbaum 72}, что при $n \geq 9$ число областей не может находиться в интервале $(3n-5, 4n-12)$. Эта гипотеза была доказана Кордовилом \cite{Cordovil 80}, Пурди \cite{Purdy 80} для $n\geq 40$ и Мартиновым \cite{Martinov 90}.
Для этого Пурди потребовалось доказать, что если для некоторого целого числа $k$
\bez
m\leq n-k \quad \text{ и} \quad n\geq 4k^2+k+1, \qquad \text{ то} \quad f\geq (k+1)(n-k).
\eez

Полное описание всех чисел отрезка $[2n-2,\frac {n(n-1)}2+1] $, которые могут реализоваться в качестве числа областей, дал Н.~Мартинов в 1993 г. \cite{Martinov 93}. В.И.~Арнольд в 2008 г., не зная о работе Мартинова, поставил задачу об описании всех возможных чисел областей "с нуля" $ $. Назвал {\it лакунами} числа интервалов $(a_i, b_i)$, где
\bez
a_i=i(n-i+1)+C_{i-1}^2, \quad b_i=(i+1)(n-i).
\eez
 С помощью неравенства
$f\geq \frac{n(n-1)}{2(m-1)}$ Арнольд \cite{Arnold 08} доказал, что число областей не может принадлежать лакуне номер $i$  для достаточно больших $n$. Однако оставалось неизвестным, содержит ли лакуна номер $i$ значения числа областей при $\frac {i^2}2 \lesssim n \lesssim i^2$ и $i \geq 3$.

\section{Теорема Мартинова о множестве возможных чисел областей}
\label{subsection theorem Martinov}

\begin{theorem} {\rm Н. Мартинов, 1993.}
\label{theorem Martinov} Нетривиальный набор из $n$ псевдопрямых на плоскости $\RP^2$ делит последнюю на $f$ областей тогда и только тогда, когда существует целое число $k$, $1\leq k \leq n-2$, такое что
\begin{equation*}
(n-k)(k+1)+C_k^2-\min\left\{n-k,C_k^2\right\}\leq f \leq (n-k)(k+1)+C_k^2.
\end{equation*}
\end{theorem}

\noindent
{\it Достаточность.} Для данного $k$, $1\leq k \leq n-2$ и числа $f$, такого что
\bez
f=(n-k)(k+1)+C_k^2\ - \ t, \qquad 0 \leq t \leq  \min\left\{n-k,C_k^2\right\}
\eez
построим набор из $n$ прямых, делящий плоскость на $f$ областей. Пусть через одну точку проходит $n-k$ прямых, остальные прямые находятся в общем положении по отношению друг к другу. При этом $t$ прямых из тех, которые проходят через одну точку, проходят через $t$ точек пересечения прямых, находящихся в общем положении. Провести так прямые возможно, т.к. $0 \leq t \leq  \min\left\{n-k,C_k^2\right\}$.
Нетрудно подсчитать, что число областей равно $f=(n-k)(k+1)+C_k^2 \ - \ t$.

\smallskip
Перед доказательством необходимости получим следствия из теоремы Мартинова.
Заметим, что объединение отрезков
\bez
\left[(n-k)(k+1)+C_k^2-\min\left\{n-k,C_k^2\right\}; \ (n-k)(k+1)+C_k^2 \right]
\eez
по $k$, $1\leq k \leq n-2$ покрывает все целые числа отрезка $[2n-2; \ C_n^2+1]$, кроме чисел интервалов $(a_i; b_i)$ --- лакун\footnote{Поскольку Арнольд не исключал тривиальные наборы, то в \cite{Arnold 08} лакуны нумеровались числами $i$, начиная с единицы. Без тривиального набора псевдопрямых лакуна $(a_1,b_1)$ пропадает, и в наших обозначениях $i \geq 2$.},
где
\bez
a_i=i(n-i+1)+C_{i-1}^2, \quad b_i=(i+1)(n-i).
\eez
Лакуна номер $i$ содержит хотя бы одно целое число тогда и только тогда, когда
\bez
b_i \geq a_i+2 \quad \lra \quad n \geq C_{i+1}^2+3.
\eez
Обозначим через $d_n$ следующее число
\bez
d_n=\max\{d \in \Z \ | \ n\geq C_{d+1}^2+3\}.
\eez
Решая квадратное неравенство, получим его явный вид
\bez
d_n=\left[\sqrt{2n-5\frac 34}- \frac 12\right].
\eez
Итак, число лакун равно $d_n-1$ и они нумеруются числами $i$, $2 \leq i \leq d_n$. По теореме Мартинова число областей не может принадлежать какой-либо лакуне. Последняя лакуна заканчивается перед числом $(n-d_n)(d_n+1)\thickapprox \sqrt{2}n^{\tfrac 32}.$ Значит, доля целых чисел отрезка $[2n-2; C_n^2+1]$, реализуемых в качестве числа областей, стремится к единице при $n \to \infty.$ Все лакуны содержатся в отрезке
$[2n-2; (n-d_n)(d_n+1)]$ и доля целых чисел этого отрезка, реализуемых в качестве числа областей, стремится к $\frac 13$ при $n \to \infty.$ Действительно, количество содержащихся в лакунах целых чисел равно
\bez
\sum_{i=2}^{d_n}(b_i-a_i-1) \ = \sum_{i=2}^{d_n} (n - C_{i+1}^2 - 2) \  \thickapprox \frac{2\sqrt{2}}3n^{\tfrac 32}.
\eez

Вставить рисунок со схематичным распределением.

\noindent
{\it Необходимость.} Напомним, что через $m$ обозначено максимальное число псевдопрямых набора, имеющих общую точку.

\bl {\rm Арнольд, \cite{Arnold 08}.}
\label{lemma Arnold}
Для нетривиальных наборов $n$ псевдопрямых
\bez
m(n-m+1) \leq f \leq m(n-m+1)+C_{n-m}^2,
\eez
причем равенство справа достигается, если помимо $m$ пересекающихся в одной точке псевдопрямых остальные $n-m$ псевдопрямых находятся в общем положении друг к другу и к коллинеарным псевдопрямым.
\el
\proof Докажем лемму индукцией по $n$, база $n=m+1$. Пусть двойное неравенство выполняется для наборов $n$ псевдопрямых. Тогда при добавлении еще одной псевдопрямой номер $n+1$ число областей увеличится на число точек пересечения добавленной псевдопрямой с предыдущими, которое не меньше $m$ и не больше $n$. Заметим, что равенство слева, вообще говоря, может не достигаться при $m< \frac n2$.
\ep

\bl
\label{lemma f v lakune}
Пусть число областей нетривиального набора $n$ псевдопрямых принадлежит лакуне номер $i,$ где $i \leq d_n$. Тогда
$m\leq i$.
\el
\proof
Предположим противное, т.е. что $m \geq i+1.$ Тогда, если $i+1 \leq m \leq n-i,$ то по лемме \ref{lemma Arnold}
\bez
f \geq m(n-m+1) \geq (i+1)(n-i)=b_i,
\eez
что противоречит тому, что число $f$ принадлежит лакуне номер $i$.
Если $m\geq n-i+1$, то $n-m \leq i-1$ и $C_{n-m}^2 \leq C_{i-1}^2$. Следовательно по лемме \ref{lemma Arnold}
\bez
f \leq m(n-m+1)+C_{n-m}^2 \leq i(n-i+1)+C_{i-1}^2 =a_i,
\eez
что противоречит тому, что число $f$ принадлежит лакуне номер $i$.
\ep

\bl {\rm Мартинов, \cite[th. 1]{Martinov 93}.}
\label{lemma Martinov}
Для нетривиальных наборов $n$ псевдопрямых и целых чисел $k$, таких что
$n\geq C_{k+1}^2+3$ и $m \leq k$ справедливо
\bez
f \geq (k+1)(n-k).
\eez
\el
Доказательство этой леммы (не оригинальное) приведено после теоремы \ref{f n m high inequality}, поскольку ее использует (само собой разумеется, теорема \ref{f n m high inequality} доказывается независимо от леммы).

Необходимость в теореме Мартинова теперь следует из лемм \ref{lemma f v lakune} и \ref{lemma Martinov}, в которой в качестве числа $k$ следует взять номер $i$ гипотетической лакуны, содержащей число областей $f$.
\ep

\section{Линейные по $t_i$ неравенства}
\label{subsection linear t_i ine}

К настоящему времени известно три линейных по $t_i$ неравенства для нетривиальных наборов $n$ прямых на проективной плоскости $\RP^2$.

{\it Неравенство Мельхиора, {\rm \cite{Melchior 40}}}. Для нетривиальных наборов псевдопрямых на плоскости $\RP^2$
\begin{equation}
\label{Melchior ine}
t_2\geq 3+\sum_{i\geq 4}(i-3)t_i.
\end{equation}

{ \it Неравенство Хирцебруха, {\rm \cite{Hirzebruch  86}}}. Для наборов комплексных прямых на комплексной проективной плоскости $\CP^2$
\begin{equation}
\label{Hirzebruch ine}
t_2+\frac 34 t_3\geq n+\sum_{i\geq 5}(2i-9)t_i \quad \text{при} \quad t_{n-1}=t_{n-2}=0.
\end{equation}

{\it Комбинаторный аналог неравенства Хирцебруха}. Для нетривиальных наборов псевдопрямых на плоскости $\RP^2$
\be
\label{combi Hirzebruch ine}
t_2+\frac 32 t_3\geq 8+\sum_{i\geq 4}\left(2i-7\frac 12\right)t_i \quad \text{при} \quad t_{n-1}=t_{n-2}=0.
\ee

Неравенство Мельхиора обращается в равенство для симплициальных наборов псевдопрямых (т.е. каждая область примыкает к трем псевдопрямым). Заметим, что неравенство Мельхиора для наборов комплексных прямых, вообще говоря, не выполняется.

Равенство в (\ref{Hirzebruch ine}) достигается, например, для набора 6 вещественных прямых, являющихся сторонами и диагоналями некоторого четырехугольника. Для такого набора $t_2=3, \ t_3=4, \ t_i=0$ для $i\geq 4$. Также неравенство (\ref{Hirzebruch ine}) обращается в равенство для набора 9 прямых с $t_4=3, \ t_3=4, \ t_2=6.$ В \cite{Hirzebruch  86} приведены примеры наборов комплексных прямых, для которых (\ref{Hirzebruch ine}) обращается в равенство.

Равенство в (\ref{combi Hirzebruch ine}) достигается для единственного с точностью до изоморфизма набора семи прямых, задаваемого  двумя точками $A$ и $B$, через каждую из которых проходит по 4 прямые набора (прямая $AB$ --- общая). Тогда $t_4=2, \ t_2=9$, $t_3=t_i=0$ при $i\geq 5$.  Для всех остальных наборов псевдопрямых имеем
\begin{equation*}
t_2+\frac 32t_3 \geq 9+\sum_{i\geq 4}\left(2i-7\frac 12\right)t_i.
\end{equation*}

\bl
\label{lemma v e f via t_i}
Для нетривиальных наборов псевдопрямых числа $v, e$ и $f$ клеточного комплекса выражаются через $t_i$:
\bez
v=\sum_{i\geq 2} t_i, \quad e=\sum_{i\geq 2}it_i, \qquad f=1+\sum_{i\geq 2} (i-1)t_i.
\eez
\proof Число $v$ выражается указанным способом по определению чисел $t_i$. Из каждой точки пересечения $i$ псевдопрямых выходит $2i$ ребер комплекса, поэтому сумма $\sum_{i\geq 2}2it_i$ суть количество ребер, посчитанных дважды. Из эйлеровой характеристики проективной плоскости имеем $v-e+f=1$, откуда следует формула для числа $f$.
\ep
\el

\begin{lemma} {\rm Мельхиор, \cite{Melchior 40}.}
\label{lemma Melchiora}
 Для нетривиальных наборов псевдопрямых
\begin{equation}
\label{iquality Melchiora}
\sum_{i\geq 2} (3-i)t_i=3+\sum_{j\geq 3}(j-3)p_j.
\end{equation}
\end{lemma}

\proof
Числа вершин $v$, ребер $e$ и областей $f$ комплекса равны
\begin{equation*}
v=\sum_{i\geq 2}t_i, \quad e=\sum_{i\geq 2}it_i=\frac 12 \sum_{j\geq 3}jp_j, \quad f=\sum_{j\geq 3}p_j.
\end{equation*}
 Тогда по формуле Эйлера для проективной плоскости получаем
\begin{equation*}
3=3f-(2e+e)+3v=3\sum_{j\geq 3}p_j-\left(\sum_{j\geq 3}jp_j+\sum_{i\geq 2}it_i\right)+3\sum_{i\geq 2}t_i=\sum_{j\geq 3}(3-j)p_j+\sum_{i\geq 2}(3-i)t_i.
\end{equation*}
\ep

Из леммы \ref{lemma Melchiora} вытекает неравенство (\ref{Melchior ine}), т.к. числа $p_j$ для $j \geq 3$ неотрицательны.

Доказательство неравенства Хирцебруха использует теорему Йо. Мияока \cite{Miyaoka 84} о числах Черна, примененную к построенному по набору комплексных прямых алгебраическому многообразию. В \cite[c. 315]{Brass} был поставлен вопрос об элементарном доказательстве неравенства Хирцебруха. Этого сделать не удалось,  неравенства Хирцебруха (\ref{Hirzebruch ine}) и его комбинаторный аналог (\ref{combi Hirzebruch ine}) не выводятся друг из друга непосредственно.
Доказательство (\ref{combi Hirzebruch ine}) опирается на следующие леммы, сформулированные для графов, ассоциированных с наборами псевдопрямых.

Нетривиальному набору из $n\geq 3$ различных псевдопрямых на вещественной проективной плоскости поставим в соответствие граф (вложенный в проективную плоскость), вершинами и ребрами которого соответственно являются точки пересечения прямых и отрезки прямых, не содержащие отличных от своих концов точек пересечения прямых. Степень любой вершины (т.е. число исходящих ребер) четна, так как точка пересечения $i$ прямых является вершиной степени $2i$. Следовательно, число вершин графа степени $2i$ равно $t_i$ для $i=2, \dots , n.$ Компоненты связности дополнения в проективной плоскости к вложенному графу суть открытые области. Любое ребро графа примыкает к двум различным областям на проективной плоскости. В графе нет петель и кратные ребра между парой вершин, если они есть, лежат на одной прямой.

\begin{lemma} {\rm (лемма о простом ребре)}
\label{lemma simple edge}
Пусть для набора $n$ прямых верно $t_n=t_{n-1}=0$. Пусть степени обоих концов некоторого ребра соответствующего графа равны четырем. Тогда из двух примыкающих к этому ребру областей хотя бы одна ограничена не менее чем четырьмя ребрами графа.
\end{lemma}

\proof Обозначим концы этого ребра через $A$ и $B$, а прямые, проходящие через точки $A$ и $B$ и отличные от прямой $AB$, через $l_1$ и $l_2$ соответственно. Обозначим точку пересечения прямых $l_1$ и $l_2$ через $C$. Предположим, что обе примыкающие к ребру $AB$ области ограничены тремя ребрами. Одно из этих ребер --- это $AB$, а два других (для каждой области) находятся на прямых $l_1$ и $l_2$ и имеют общую точку. Точка $C$ --- это единственная общая точка прямых $l_1$ и $l_2$, поэтому каждая из двух примыкающих к ребру $AB$ областей ограничена ребром с концами в точках $A$ и $C$, ребром с концами в точках $B$ и $C$ и самим ребром $AB$. Следовательно, на прямой $l_1$ есть ровно две точки пересечения с остальными прямыми набора, и это точки $A$ и $C$. На проективной плоскости любые две различные прямые пересекаются в единственной точке, поэтому каждая прямая из набора, кроме $l_1$, проходит или через точку $A$, или через точку $C$. Степень точки $A$ равна четырем, т.е. через точку $A$ проходит, не считая $l_1$, только прямая $AB$. Тогда через точку $C$ вместе с прямой $l_1$ проходит  $n-1$ прямых, что противоречит условию $t_{n-1}=0$. Следовательно, обе примыкающие к ребру $AB$ области не могут быть ограничены тремя ребрами графа каждая. $\quad \square$

\begin{lemma} {\rm (лемма об оценке $t_2$ сверху).}
\label{lemma t_i p_j inequality}
 При условии $t_n=t_{n-1}=t_{n-2}=0$ справедливо неравенство
\begin{equation}
\label{t_i p_j inequality}
2t_2\leq 1+ 3p_4+\sum_{j\geq 5}jp_j+\sum_{i\geq 3}\left(i-\frac 32\right)t_i.
\end{equation}
\end{lemma}

\proof Для соответствующего набору прямых графа обозначим через $x$ число ребер, оба конца которых имеют степень 4, а через $y$ --- число ребер, оба конца которых имеют степень не менее 6.

{\it Шаг 1.} Всего в графе $\sum_{i\geq 2}it_i$ ребер, поэтому число ребер, один конец которых имеет степень 4, а другой не меньшую чем 6, равно $\sum_{i\geq 2}(it_i) \ -x-y$. Каждая вершина графа степени 4 является концом четырех ребер, хотя бы один конец которых имеет степень 4. Поэтому суммарное по всем ребрам число их концов степени 4 равно $4t_2$ и равно
\begin{equation}
\label{x y t_i inequality}
4t_2=2x+\sum_{i\geq 2}it_i \ - \ x-y \qquad \Longrightarrow \qquad x=2t_2+y-\sum_{i\geq 3}it_i.
\end{equation}

{\it Шаг 2.} Предположим, что существуют две различные точки $A$ и $B$, такие что любая прямая из набора проходит через хотя бы одну из них. Обозначим через $a$ и $b$ число прямых набора, проходящих через точки $A$ и $B$ соответственно. Возможны два случая.

\noindent
(i) Прямая $AB$ не принадлежит набору $n$ прямых. Тогда $a+b=n$  и из условия $t_{n-2}=0$ следует, что $a\geq 3$ и $b \geq 3$. В этом случае
\begin{equation*}
t_2=ab, \qquad \sum_{i\geq 3}\left(i-\frac 32\right)t_i=a+b-3, \qquad p_4=ab-a-b+3, \qquad p_i=0
 \quad \text{при} \quad i\geq 5.
\end{equation*}
Теперь неравенство (\ref{t_i p_j inequality}) проверяется непосредственно:
\begin{equation*}
2ab\leq 1+3(ab-a-b+3)+a+b-3 \qquad \quad \Leftrightarrow \qquad \quad (a-2)(b-2)+3\geq 0.
\end{equation*}
(ii) Прямая $AB$ принадлежит набору $n$ прямых. Тогда $a+b=n+1$ и из условия $t_{n-2}=0$ следует, что $a\geq 4$ и $b \geq 4$. В этом случае
\begin{equation*}
t_2=ab-a-b+1, \quad \sum_{i\geq 3} \left(i-\frac 32\right)t_i=a+b-3, \quad p_4=ab-2a-2b+4, \qquad p_i=0
 \quad \text{при} \quad i\geq 5.
\end{equation*}
Теперь неравенство (\ref{t_i p_j inequality}) проверяется непосредственно:
\begin{equation*}
2(ab-a-b+1)\leq 1+3(ab-2a-2b+4)+a+b-3 \qquad \Leftrightarrow \qquad (a-3)(b-3)-1\geq 0.
\end{equation*}
В дальнейшем доказательстве леммы (а именно, в шагах 5 и 6) будем считать, что не существует двух различных точек, таких что любая прямая набора проходит через хотя бы одну из них.

{\it Шаг 3.} Для данного графа рассмотрим множество $F$ областей (т.е. компонент связности дополнения к прямым), каждая из которых ограничена не менее чем четырьмя ребрами и граница которой содержит хотя бы одну вершину степени 4 (внутри областей точек графа нет). Для области $\Gamma \in F$ обозначим через $x(\Gamma)$ число ограничивающих $\Gamma$ ребер, оба конца которых имеют степень 4. Для области $\Gamma \in F$  обозначим через $s(\Gamma)$ число ее вершин (т.е. вершин на границе $\Gamma$) степени не менее чем 6. Положим
\begin{equation*}
\delta(\Gamma) = \left\{
   \begin{array}{ll}
     0,  & \hbox{если \ \ $s(\Gamma)\geq 1$;} \\
     1,  & \hbox{если \ \ $s(\Gamma)=0$.}
   \end{array}
 \right.
\end{equation*}
Докажем, что если область $\Gamma$ ограничена $j$ ребрами, то
\begin{equation}
\label{s j x inequality}
s(\Gamma)\leq (j-1)-x(\Gamma)+\delta(\Gamma). \end{equation}
Рассмотрим три случая.

\noindent
(i) $x(\Gamma)=0$. Тогда $s(\Gamma)\leq j-1$, так как на границе $\Gamma$ есть вершина степени 4.

\noindent
(ii) $x(\Gamma)=j$. Тогда $s(\Gamma)=0$ и $\delta(\Gamma)=1$.

\noindent
(iii) $0<x(\Gamma)<j$. Рассмотрим отдельно границу $\Gamma$, состоящую из $j$ ребер. Среди них выберем $x(\Gamma)$ ребер, оба конца которых имеют степень 4. Пусть эти $x(\Gamma)$ ребер образуют на границе $\Gamma$ ровно $z(\Gamma)$ компонент связности, каждая компонента --- это несколько подряд идущих выбранных ребер. Из $x(\Gamma)>0$ следует, что $z(\Gamma)\geq 1$, а из $x(\Gamma)<j$ следует, что каждая компонента не замкнута (т.е. гомеоморфна отрезку). В каждой компоненте число вершин степени 4 на единицу больше числа ребер этой компоненты. Все вершины и ребра различных компонент связности различны, поэтому граница области $\Gamma$ содержит не менее $x(\Gamma)+z(\Gamma)$ вершин степени 4. Так как $z(\Gamma)\geq 1$, то $s(\Gamma)\leq j-1-x(\Gamma).$

Обозначим через $s$ сумму $s=\sum_{\Gamma \in F}s(\Gamma)$. Суммируя (\ref{s j x inequality}) по всем областям $\Gamma \in F$, получим
\begin{equation}
\label{s p_j x inequality}
s\leq \sum_{j\geq 4}(j-1)p_j-\sum_{\Gamma \in F}x(\Gamma)+\sum_{\Gamma \in F}\delta(\Gamma).
\end{equation}

 {\it Шаг 4.} Покрасим в красный цвет ребра графа, оба конца которых имеют степень 4 и обе примыкающие области к которым ограничены не менее чем четырьмя ребрами каждая (т.е. обе примыкающие области из $F$). Обозначим число красных ребер через $a$. Тогда по лемме \ref{lemma simple edge} (о простом ребре) число $x-a$ равно числу ребер, оба конца которых имеют степень 4 и ровно одна из двух примыкающих областей ограничена не менее чем четырьмя ребрами (т.е. одна примыкающая область из $F$).
 Следовательно,
\begin{equation}
\label{x(G) equality}
\sum_{\Gamma \in F}x(\Gamma)=x+a.
\end{equation}

Покрасим в синий цвет четырехугольные области, все вершины которых имеют степень 4. Докажем, что к каждой синей области примыкает не менее двух красных ребер. Для этого выведем из условия $t_{n-2}=0$ аналогично лемме \ref{lemma simple edge} (о простом ребре), что в каждой паре противоположных ребер любой синей области есть хотя бы одно красное ребро. Предположим противное, что оба ребра $AB$ и $CD$ некоторой синей области $ABCD$ не красные. Тогда к ребрам $AB$ и $CD$ примыкают треугольные области $ABH$ и $CDG$, причем точка пересечения прямых $BC$ и $AD$ совпадает и с точкой $H$, и с точкой $G$. Следовательно, через эту точку $G=H$ проходят $n-2$ прямые набора (все прямые кроме $AB$ и $CD$), что противоречит условию $t_{n-2}=0$.

Обозначим число синих областей через $p$.
Сумма $\sum_{\Gamma \in F} \delta(\Gamma)$ равна количеству областей, каждая из которых ограничена не менее чем четырьмя ребрами и все вершины которой имеют степень 4. Поэтому
\begin{equation}
\label{d p p_j ine}
\sum_{\Gamma \in F}\delta(\Gamma) \leq p+\sum_{j\geq 5}p_j.
\end{equation}

Обозначим через $\varphi$ число пар $(C,\kappa)$ синих областей $C$ и красных ребер $\kappa$ на границе области $C$.
Так как к любой из $p$ синих областей примыкает не менее двух красных ребер, то $\varphi \geq 2p$. С другой стороны, каждое красное ребро примыкает к не более двум синим областям и поэтому  $2a\geq \varphi.$ Следовательно, $a\geq p$.
Итак, из (\ref{s p_j x inequality}), (\ref{x(G) equality}), (\ref{d p p_j ine}) и неравенства $a\geq p$ следует, что
\begin{equation}
\label{s p_4 x p_j ine}
s\leq 3p_4-x+\sum_{j\geq 5}jp_j.
\end{equation}

{\it Шаг 5.} Рассмотрим произвольную вершину $V$ степени не менее $ 6$ и удалим из графа все ребра, лежащие на проходящих через точку $V$ прямых (соответственно изменятся степени оставшихся вершин, а некоторые вершины, возможно, исчезнут). Обозначим полученный граф через $G(V)$, а исходный граф через $G$. После шага 2 достаточно рассматривать только те наборы прямых, для которых не существует двух точек, таких что любая прямая набора проходит через хотя бы одну из них. Для таких наборов прямых граф $G(V)$ имеет хотя бы две различные вершины и каждая область проективной плоскости, образованная графом $G(V)$, ограничена не менее чем тремя ребрами графа $G(V)$. Значит, точка $V$ находится внутри некоторой области, образованной графом $G(V)$, граница которой есть $d$-угольник с вершинами $A_1, \dots , A_d$ при $d\geq 3$ (занумерованными в порядке следования).
Докажем, что для любой вершины $V\in G$ имеет место хотя бы одно из следующих утверждений.

\noindent
(а) Вершина $V$ соединена ребрами графа $G$ с не менее чем тремя вершинами из множества $\{ A_1, \dots , A_d\}$.

\noindent
(б) Вершина $V$ соединена ребрами графа $G$ с двумя вершинами из множества $\{A_1, \dots, A_d\}$ и является вершиной границы некоторой области из множества $F$.

\noindent
(в) Вершина $V$ является вершиной границ не менее двух областей из множества $F$.

 Пусть для какого-то $i$, $1\leq i \leq d$ отрезок $VA_i$ не является ребром графа $G$. Тогда интервал $VA_i$ принадлежит образованной графом $G$ области, граница которой содержит вершины $V$ и $A_i$, и состоит из не менее четырех ребер графа $G$. Если эта область не из множества $F$, то ее граница есть многоугольник с вершинами $V, A_k, A_{k+1}, \dots, A_i, \dots, A_l$ для некоторых чисел $k$ и $l$, т.е. все вершины границы, кроме $V$, суть точки $A_k, \dots, A_l$ и отрезки $VA_k$ и $VA_l$ являются ребрами графа $G$.

Предположим, что для вершины $V$ утверждение (а) не выполняется. Тогда возможны два случая.

\noindent
(i) Вершина $V$ соединена ребрами графа $G$ ровно с двумя вершинами из множества $\{A_1, \dots,A_d\}$.

\noindent
(ii) Вершина $V$ соединена ребрами графа $G$ с не более одной вершиной из множества $\{A_1, \dots,A_d\}$.

В случае (i) обозначим  эти  две  вершины  через $A_k$ и $A_l$. Тогда среди остальных $d-2~\geq~1$ вершин множества $\{A_1, \dots,A_d\}$ найдется вершина $A_q$, такая что точка $V$ содержится в треугольнике $A_qA_kA_l$ (а треугольник лежит в замыкании области $A_1,\dots, A_n$). Тогда отрезок $VA_q$ не является ребром графа $G$ и интервал $VA_q$ содержится в области из множества $F$. Следовательно, в случае (i) выполняется утверждение (б).

Cлучай (ii). Предположим, что вершина $V$ соединена ребром графа $G$ с некоторой вершиной $A_k$ из множества $\{A_1, \dots,A_d\}$. Тогда прямая $VA_k$ разбивает множество вершин $\{A_1, \dots,A_d\}\setminus\{A_k\}$ на два непустых подмножества, находящихся по разные стороны (в замыкании области $A_1 \dots A_d$) от прямой $VA_k$. Тогда по обе стороны от прямой $VA_k$ найдется область из $F$, граница которой содержит точку $V$. Следовательно, будет выполнено утверждение (в).
Если вершина $V$ не соединена ребром графа $G$ ни с какой вершиной из множества $\{A_1, \dots,A_d\}$, то вместо прямой $VA_k$ возьмем любую прямую набора, проходящую через точку $V$. Аналогично получим, что для вершины $V$ выполнено утверждение (в).

{\it Шаг 6.} Для произвольной вершины $V$ (исходного) графа степени не меньше 6 обозначим через $s'(V)$ сумму числа подходящих к вершине $V$ ребер, другой конец которых имеет степень не меньше 6, и удвоенного числа примыкающих к $V$ областей из множества $F$. Так как для любой вершины $V$ степени не меньше чем 6 выполняется хотя бы одно из утверждений (а)---(в) предыдущего шага, то $s'(V) \geq 3$. Обозначим через $s'$ сумму \begin{equation*}
s'=\sum_{\deg(V)\geq 6}s'(V)
\end{equation*}
чисел $s'(V)$ для всех вершин $V$ степени не меньше чем 6. Так как $s'(V) \geq 3$, то
$s' \geq 3\sum_{i\geq 3}t_i$. С другой стороны, $s'=2y+2s$. Следовательно,
\begin{equation}
\label{y+s ine}
y+s\geq \frac 32\sum_{i\geq 3}t_i.
\end{equation}
Из (\ref{x y t_i inequality}) и
(\ref{s p_4 x p_j ine}) получаем
\begin{equation*}
3p_4+\sum_{j\geq 5}jp_j \ - \ s \geq  x=2t_2+y-\sum_{i\geq 3}it_i \quad \Longrightarrow \quad 3p_4+\sum_{j\geq 5}jp_j +\sum_{i\geq 3}it_i\ - \ 2t_2\geq y+s.
\end{equation*}
Из последнего неравенства и (\ref{y+s ine}) следует, что
\begin{equation*}
3p_4+\sum_{j\geq 5}jp_j+\sum_{i\geq 3}it_i \ - \ 2t_2 \geq \frac 32 \sum _{i\geq 3}t_i \quad \Longrightarrow \quad 3p_4+\sum_{j\geq 5}jp_j+\sum_{i\geq 3}\left(i-\frac 32\right)t_i\geq 2t_2. \quad \square
\end{equation*}

\begin{theorem}
\label{theor combi Hirzebruch} Для нетривиальных наборов псевдопрямых с $t_{n-1}=t_{n-2}=0$
\bez
t_2+\frac 32 t_3\geq 8+\sum_{i\geq 4}\left(2i-7\frac 12\right)t_i.
\eez
\end{theorem}

\proof По лемме Мельхиора \ref{lemma Melchiora}  имеем
\begin{equation*}
\sum_{i\geq 2}(9-3i)t_i=9+ 3p_4+\sum_{j\geq 5}(3j-9)p_j.
\end{equation*}
По лемме \ref{lemma t_i p_j inequality} получаем
\begin{equation*}
3p_4+\sum_{j\geq 5}jp_j\geq 2t_2-\sum_{i\geq 3}\left(i-\frac 32\right)t_i\ -\ 1.
\end{equation*}
Заметим, что $3j-9\geq j$ при $j\geq 5.$ Следовательно, т.к. $p_j\geq 0$, справедливо \begin{equation*}
\sum_{i\geq 2}(9-3i)t_i\geq 9+2t_2-\sum_{i\geq 3}\left(i-\frac 32\right)t_i\ -\ 1 \quad \Leftrightarrow \quad t_2+\frac 32t_3 \geq 8+\sum_{i\geq 4}\left(2i-7\frac 12\right)t_i. \quad \square
\end{equation*}

\begin{remark} {\rm (идея А.Т.Фоменко.)} Перенесем все члены неравенств Мельхиора (\ref{Melchior ine}), Хирцебруха (\ref{Hirzebruch ine}) и (\ref{combi Hirzebruch ine}) в большую часть и составим из получившихся коэффициентов при $t_2,\dots, t_n$ три вектора в $\mathbb{R}^{n-1}$. А именно, координата номер $i=1,\dots, n-1$ векторов равна коэффициенту при $t_{i+1}$ в соответствующем неравенстве:

\begin{equation*}
\overrightarrow{N_1}=\left(1,0,-1,\dots,3-n\right), \ \overrightarrow{N_2}=\left(1, \frac 34, 0, -1, \dots, 9-2n\right), \ \overrightarrow{N_3}=\left(1, \frac 32, -\frac 12, \dots, 7\frac 12 - 2n\right). \end{equation*}
Найдем асимптотику длин векторов $\overrightarrow{N_1},\  \overrightarrow{N_2}, \  \overrightarrow{N_3} $ при $n\rightarrow \infty:$
\begin{equation*}
|\overrightarrow{N_1}|=\frac 1{\sqrt {3}}n^{\tfrac 32}\left(1+O\left(\frac 1n\right)\right), \ |\overrightarrow{N_2}|=\frac 2{\sqrt {3}}n^{\tfrac 32}\left(1+O\left(\frac 1n\right)\right), \ |\overrightarrow{N_3}|=\frac 2{\sqrt {3}}n^{\tfrac 32}\left(1+O\left(\frac 1n\right)\right).
\end{equation*}
Докажем, что углы между векторами $\overrightarrow{N_1}$ и $\overrightarrow{N_2}$, а также между векторами $\overrightarrow{N_2}$ и $\overrightarrow{N_3}$ стремятся к нулю при $n\rightarrow \infty.$ Так как
\begin{equation*}
|\overrightarrow{N_1} - \overrightarrow{N_2}|=\frac 1{\sqrt {3}}n^{\tfrac 32}\left(1+O\left(\frac 1n\right)\right)\quad \text{и} \quad |\overrightarrow{N_2} - \overrightarrow{N_3}|=\frac 32 n^{\tfrac 12}\left(1+O\left(\frac 1n\right)\right) \quad \text{при} \quad n\rightarrow \infty,
\end{equation*}
то по теореме косинусов для треугольника со сторонами $\overrightarrow{N_1}$, $\overrightarrow{N_2}$ и треугольника со сторонами $\overrightarrow{N_2}$, $\overrightarrow{N_3}$ получаем:
\begin{equation*}
\cos \ \angle\left(\overrightarrow{N_1},\overrightarrow{N_2}\right)=
 1- O\left(\frac 1n\right)\quad \text{и} \quad \cos \ \angle\left(\overrightarrow{N_2},\overrightarrow{N_3}\right)=
  1- O\left(\frac 1n\right) \quad \text{при} \quad  n \rightarrow \infty.
\end{equation*}
\end{remark}

\section{Нижние оценки числа областей}
\label{subsection ocenki}

\begin{theorem}
\label{f n m high inequality}
Пусть для нетривиального набора $n$ различных псевдопрямых на вещественной проективной плоскости максимальное число псевдопрямых, имеющих общую точку, равно $m$. Пусть $T \geq m$. Тогда
\begin{gather}
\label{f n m medium inequlity}
f \geq 2\frac{n^2-n+2T}{T+3}, \\
 \quad f\geq
\frac{\left(3m-8\tfrac 12\right)(n^2-n)+(9m^2-21m+1)}{m^2+3m-15} \quad \text{при } \quad 12 \leq m<n-2. \label{f n m strong ine b}
\end{gather}
Если для нетривиального набора $n$ различных прямых на плоскости $\RP^2$ максимальное число псевдопрямых, имеющих общую точку, равно $m$. Тогда
\be
f\geq \frac{(3m-10)n^2+(m^2-6m+12)n}{m^2+3m-18} \ +\ 1 \quad \text{при} \quad 5 \leq m<n-2 , \label{f n m strong ine a}
\ee
\end{theorem}

Заметим, что для наборов прямых выполняются все три неравенства теоремы. Также заметим, что (\ref{f n m medium inequlity}) следует из (\ref{f n m strong ine b}) при $6 \leq m<n-2$.

\smallskip
 Автором был предложен метод получения оценок числа областей $f$ снизу с помощью линейных неравенств на числа $t_i$. Впервые метод был применен в \cite{Shnurnikov 10} к неравенству (\ref{Melchior ine}), что дало оценку (\ref{f n m medium inequlity}). Неравенства (\ref{f n m strong ine a}) и (\ref{f n m strong ine b}) получаются из неравенств (\ref{Hirzebruch ine}) и (\ref{combi Hirzebruch ine}) соответственно.

 \smallskip
Зафиксируем число прямых $n$ и максимальное число $m$ прямых, имеющих общую точку. По лемме \ref{lemma v e f via t_i}
\begin{equation}
\label{f via t_i in a proof}
f-1=\sum_{i=2}^{m}(i-1)t_i.
\end{equation}
Получим формулу
\begin{equation}
\label{n, m, t_i equality}
n(n-1)=\sum_{i= 2}^{m}i(i-1)t_i
\end{equation}
подсчетом числа $\frac{n(n-1)}2$ пар прямых. Каждой точке пересечения $i$ прямых поставим в соответствие $\frac{i(i-1)}2$ пар прямых, проходящих через эту точку. Поскольку любые две прямые пересекаются в одной точке, то в сумме $\sum_{i\geq 2}\frac{i(i-1)}2t_i$ каждая пара прямых учитывается ровно один раз.

Рассмотрим линейное по $t_i$ неравенство вида
\begin{equation}
\label{alpha_i t_i general}
\sum_{i\geq 2}\alpha_it_i\geq \alpha_0, \end{equation}
где $\alpha_0, \alpha_2, \alpha_3, \dots, \alpha_n$ --- константы при фиксированном $n$.
Например, для неравенства (\ref{Hirzebruch ine}) эти константы равны
\begin{equation*}
\alpha_0=n,\quad \alpha_2=1, \quad \alpha_3=\frac 34, \quad \alpha_4=0, \quad \alpha_i=9-2i \quad \text{при }\quad i\geq 5. \end{equation*}

Подберем такие положительные коэффициенты $c_1$ и $c_2$ (постоянные при фиксированном $m$), что
\begin{equation}
\label{c_1 c_2 ine}
c_1i(i-1)+c_2\alpha_i\leq i-1 \qquad \text{для всех }\qquad 2\leq i \leq m.
\end{equation}
Умножим неравенство (\ref{c_1 c_2 ine}) для каждого $i$ на $t_i$ и просуммируем по $i=2, \dots, m$. Поскольку числа $t_i$ неотрицательны, то получится неравенство \begin{equation*}
c_1\sum_{i=2}^{m}i(i-1)t_i\ +\  c_2\sum_{i=2}^{m}\alpha_it_i \leq \sum_{i=2}^{m}(i-1)t_i \quad \Leftrightarrow \quad c_1n(n-1)+c_2\sum_{i=2}^{m}\alpha_it_i \leq f-1.
\end{equation*}
Так как $c_2>0$, то из последнего неравенства, из (\ref{alpha_i t_i general}) и того, что $t_k=0$ при $k>m$, следует
\begin{equation}
\label{f c_1 c_2 ine}
f\geq c_1n(n-1)+c_2\alpha_0+1. \end{equation}
для положительных констант $c_1$ и $c_2$, удовлетворяющих системе неравенств (\ref{c_1 c_2 ine}).

Теперь применим метод к известным линейным по $t_i$ неравенствам.

\prove{теоремы \ref{f n m high inequality}}
Докажем первое неравенство. Запишем неравенство Мельхиора (\ref{Melchior ine}) в виде (\ref{alpha_i t_i general}) с коэффициентами
\begin{equation*}
 \alpha_0=3,\qquad \alpha_i=3-i \quad \text{при }\quad i\geq 2.
 \end{equation*}
Введем положительные множители
\bez
c_1=\frac 2{T+3}\quad \text{и }\quad c_2=\frac{T-1}{T+3}.
\eez
Рассмотрим квадратный многочлен относительно $i$: $w(i)=c_1 i(i-1)+c_2(3-i)-(i-1)$. Так как значения многочлена $w(i)$ на концах отрезка $[2,T]$ равны нулю и старший коэффициент $c_1>0,$ то $w(i)\leq 0$ для всех  $2\leq i \leq T.$ Отсюда и из (\ref{f c_1 c_2 ine}) получаем
\bez
f \geq 2\left(\frac{n^2-n+2T}{T+3}\right).
\eez

Докажем второе неравенство теоремы. Так как $m<n-2$, то $t_{n-1}=t_{n-2}=0$ и выполняется неравенство
 (\ref{combi Hirzebruch ine}). Запишем (\ref{combi Hirzebruch ine}) в виде (\ref{alpha_i t_i general}) с коэффициентами \begin{equation*}
 \alpha_0=8,\quad \alpha_2=1, \quad \alpha_3=\frac 32, \quad \alpha_i=7\frac 12-2i \quad \text{при }\quad i\geq 4. \end{equation*}
Так как $m\geq 12$, то система неравенств (\ref{c_1 c_2 ine}) принимает вид
 \begin{equation}
 \label{c_1 c_2 b ine}
 1\geq 2c_1+c_2, \quad 2\geq 6c_1+\frac 32 c_2, \qquad  i-1\geq c_1i(i-1)-c_2\left(2i-7\frac 12\right) \quad \text{для} \quad 4\leq i \leq m. \end{equation}
Докажем, что для констант (положительных при $m\geq 12$)
\begin{equation*}
c_1=\frac{3m-8\tfrac 12}{m^2+3m-15}, \qquad c_2=\frac{m^2-3m+2}{m^2+3m-15}
\end{equation*}
выполняются неравенства (\ref{c_1 c_2 b ine}). Во-первых,  $2c_1+c_2=1$. Во-вторых,
\begin{equation*}
6c_1+\frac 32c_2=3c_1+\frac 32(2c_1+c_2)=3c_1+\frac 32 \leq 2 \quad \Leftrightarrow \quad c_1\leq \frac 16 \quad \Leftrightarrow \quad (m-12)(m-3)\geq 0. \end{equation*}
В-третьих, рассмотрим квадратный трехчлен \begin{equation*}
R(i)=c_1i(i-1)-c_2\left(2i-7\frac 12\right)-(i-1).
\end{equation*}
Разложим трехчлен $R(i)$ на множители
\begin{equation*}
R(i)=\frac{(i-m)\left(\left(3m-8\tfrac12\right)i- \left(8\tfrac 12m-19\tfrac 12\right)\right)}{m^2+3m-15}.
\end{equation*}
Заметим, что $R(i)\leq 0$ при $4\leq i \leq m$, т.к. для $i\geq 4$ выполняется \begin{equation*}
\left(3m-8\frac 12\right)i-\left(8\frac 12m-19\frac 12\right)\geq 3\frac 12m-14\frac 12>0 \quad \text{при }\quad m\geq 12. \end{equation*}
Следовательно, для выбранных констант $c_1$ и $c_2$ справедлива система неравенств (\ref{c_1 c_2 b ine}) и неравенство (\ref{f c_1 c_2 ine}) при $\alpha_0=8$ имеет вид
\begin{equation*}
f\geq \frac{\left(3m-8\tfrac 12\right)(n^2-n)+(9m^2-21m+1)}{m^2+3m-15} \qquad \text{при }\qquad 12 \leq m < n-2.
\end{equation*}

Докажем третье неравенство.
Так как $m<n-2$, то $t_n=t_{n-1}=t_{n-2}=0$ и выполняется неравенство (\ref{Hirzebruch ine}). Запишем (\ref{Hirzebruch ine}) в виде (\ref{alpha_i t_i general}) с коэффициентами \begin{equation*}
\alpha_0=n,\quad \alpha_2=1, \quad \alpha_3=\frac 34, \quad \alpha_4=0, \quad \alpha_i=9-2i \quad \text{при }\quad i\geq 5. \end{equation*}
Так как $m\geq 5$, то система неравенств (\ref{c_1 c_2 ine}) принимает вид
 \begin{equation}
 \label{c_1 c_2 a ine}
 1\geq 2c_1+c_2, \quad 2\geq 6c_1+\frac 34 c_2, \quad 3\geq 12c_1, \quad
i-1\geq c_1i(i-1)-c_2(2i-9) \ \text{для} \  5\leq i \leq m.
\end{equation}
Докажем, что для констант (положительных при $m\geq 5$)
\begin{equation*}
 c_1=\frac{3m-10}{m^2+3m-18}, \qquad c_2=\frac{m^2-3m+2}{m^2+3m-18}
\end{equation*}
выполняются неравенства
(\ref{c_1 c_2 a ine}). Во-первых, $2c_1+c_2=1$. Во-вторых,
\begin{equation*}
c_1\leq \frac 14  \quad \Leftrightarrow \quad (m-4)(m-5)+2\geq 0.
\end{equation*}
В-третьих,
\begin{equation*}
6c_1+\frac 34c_2=\frac 34(2c_1+c_2)+4\frac 12c_1\leq \frac 34+\frac 98<2. \end{equation*}
Наконец, рассмотрим квадратный трехчлен \begin{equation*}
Q(i)=c_1i(i-1)-c_2(2i-9)-(i-1). \end{equation*}
Разложим трехчлен $Q(i)$ на множители
\begin{equation*}
Q(i)=\frac{(i-m)\left((3m-10)i-(10m-24)\right)}{(m-3)(m+6)}. \end{equation*}
Заметим, что $Q(i)\leq 0$ при $5\leq i \leq m$, так как, если $m\geq 6$, то \begin{equation*}
(3m-10)i-(10m-24)\geq 5m-26>0,
\end{equation*}
а если $m=5$, то тогда $i=5$ и $Q(5)=0$. Итак, для выбранных $c_1$ и $c_2$ справедлива система неравенств (\ref{c_1 c_2 a ine}) и поэтому из (\ref{f c_1 c_2 ine}) при $\alpha_0=n$ следует
\begin{equation*}
f\geq \frac{(3m-10)n^2+(m^2-6m+12)n}{m^2+3m-18} \ +\ 1 \quad \text{при}\quad 5 \leq m < n-2. \end{equation*}
\begin{flushright}
$\square$
\end{flushright}

\paragraf{Завершение доказательства теоремы Мартинова}

В качестве применений полученных неравенств, докажем лемму \ref{lemma Martinov}:

{\it
Для нетривиальных наборов $n$ псевдопрямых и целых чисел $k$, таких что
$n\geq C_{k+1}^2+3$ и $m \leq k$, справедливо
\be
f \geq (k+1)(n-k).
\ee
}
Для наборов прямых эту лемму можно алгебраически вывести из первого и третьего неравенств теоремы \ref{f n m high inequality}. Для наборов псевдопрямых нам понадобится дополнительное построение.
Рассмотрим четыре случая.

(1)   $m<k$,

(2) $m=k$ и $2\leq k \leq 5$,

(3) $m=k\geq 6$ и $t_k=1$,

(4) $m=k\geq 6$ и $t_k\geq 2.$

{\it Случаи 1 --- 3.} В случае 1 из набора псевдопрямых удаляется любая псевдопрямая, в случае 3 --- псевдопрямая, проходящая через точку пересечения $k$ псевдопрямых. Число областей при этом только уменьшится. Применим первое неравенство теоремы \ref{f n m high inequality} для полученных семейств из $n-1$ псевдопрямых с максимальной кратностью не более $k-1$ в случаях 1 и 3 и для исходного набора из $n$ псевдопрямых с максимальной кратностью $k$ в случае 2.
\bez
\text{Случаи 1,3:} \quad f\geq 2\frac{(n-1)^2-(n-1)+2(k-1)}{k+2}; \qquad \qquad \text{случай 2:} \quad f\geq 2\frac {n^2-n+2k}{k+3}.
\eez
Докажем следующие неравенства при условии $n\geq \frac{k^2+k}2+3:$
\bez
2\frac{(n-1)^2-(n-1)+2(k-1)}{k+2}\geq (k+1)(n-k), \qquad \quad 2\frac {n^2-n+2k}{k+3} \geq (k+1)(n-k) \quad \text{при} \quad k\leq 5.
\eez
Эти неравенства равносильны неравенствам
\bez
q(n)=n^2-n\frac{k^2+3k+8}2+\frac{k^3+3k^2+6k}2 \geq 0 \quad \text{и} \quad s(n)=n^2-n\frac{k^2+4k+5}2+\frac{k^3+4k^2+7k}2 \geq 0.
\eez
Левые части двух последних неравенств суть квадратные трехчлены $q(n)$ и $s(n)$ относительно $n$, для проверки неотрицательности которых при $n\geq \frac{k^2+k}2+3$ достаточно установить неотрицательность значений $q(\frac{k^2+k}2+3)$ и $s(\frac{k^2+k}2+3)$:
\bez
q\left(\frac{k^2+k}2+3\right)=\frac{(k-3)(k+2)}2\geq 0 \quad \text{и} \quad s\left(\frac{k^2+k}2+3\right)=\frac{(6-k)(k^2+1)-2k}4.
\eez
Заметим, что в случаях 1 и 3 верно $k\geq 3,$ поэтому $q(\frac{k^2+k}2+3)\geq 0.$ При $k\leq  5$ верно $6-k\geq 1$, и поэтому $s(\frac{k^2+k}2+3)\geq 0.$ Итак, во всех случаях 1 --- 3  получили неравенство $f\geq (k+1)(n-k).$

{\it Случай 4.} Дано $m=k\geq 6$ и $t_k\geq 2.$ Тогда
найдутся хотя бы две точки, в каждой из которых пересекается $k$
псевдопрямых. Обозначим эти точки через $P$ и $Q$. Рассмотрим
последовательность семейств псевдопрямых $A_0,A_1,\ldots,A_i,$ в которой
каждое следующее семейство получается из предыдущего добавлением
одной псевдопрямой из исходного набора $A$. Возможны два случая, в зависимости от
которых мы определим $A_0$ и $A_i$.

(а) Прямая $PQ$ не принадлежит семейству $A$. В качестве $A_0$ возьмем $2k$-элементное множество псевдопрямых, проходящих через точки $P$ и $Q$, а в качестве $A_i$ --- семейство $A$. Тогда получается $I=n-2k$.

(б) Псевдопрямая $PQ$ принадлежит семейству $A$. За конфигурацию $A_0$ возьмем $2k-2$ псевдопрямые, проходящие через точки $P$ и $Q$ и отличные от псевдопрямой $PQ$. Пусть семейство $A_i$ образуют все псевдопрямые из $A$, кроме псевдопрямой $PQ$. В этом случае $I=n-2k+1$.

Будем рассматривать оба случая одновременно, пока не получим неравенство на $f(A_i)$. Введем $R(A_i)=2e(A_i)-3f(A_i)$ --- разницу между удвоенным числом ребер и утроенным числом областей конфигурации $A_i$. Для всех $i$ верно $R(A_i)\geq 0$, так как все области ограничены хотя бы тремя дугами. Пусть семейство $A_i$ получается из $A_{i-1}$ добавлением псевдопрямой $l_i$. Пусть псевдопрямая $l_i$ проходит через $v_i$ точек пересечения псевдопрямых конфигурации $A_{i-1}$ и еще пересекает $u_i$ псевдопрямых из $A_{i-1}$ в $u_i$ точках. Тогда в конфигурации $A_i$ на псевдопрямой $l_i$ лежит $u_i$ точек пересечения кратности два и $v_i$ точек пересечения кратности хотя бы три. Следовательно, псевдопрямая $l_i$ состоит (в конфигурации $A_i$) из $u_i+v_i$ дуг, каждая из которых делит некоторую область, высекаемую семейством $A_{i-1}$, на две части. Поэтому $f(A_i)-f(A_{i-1})=u_i+v_i$. Кроме того, псевдопрямая $l_i$ делит надвое $u_i$ ребер из конфигурации $A_{i-1}$, поэтому $e(A_i)-e(A_{i-1})=2u_i+v_i$. Из двух последних формул следует неравенство
\bez
R(t_i)-R(A_{i-1})=u_i-v_i.
\eez
Пусть псевдопрямая $l_i$ проходит через $w_i$ точек пересечения псевдопрямых семейства $A_0.$ Тогда псевдопрямая $l_i$ пересекает псевдопрямые семейства $A_0$ в $n(A_0)-w_i$ точках. Пусть $z_i$ и $x_i$ обозначают количества не лежащих на псевдопрямых семейства $A_0$ точек пересечения псевдопрямой $l_i$ с псевдопрямыми семейства $A_{i-1}$, которые принадлежат ровно одной и хотя бы двум псевдопрямым семейства $A_{i-1}$ соответственно. Все оставшиеся точки пересечения псевдопрямой $l_i$ с ровно одной псевдопрямой из $A_{i-1}$ лежат на прямых из $A_0$, поэтому  оставшихся точек не более $n(A_0)-2w_i$. Следовательно, $u_i\leq n(A_0)-2w_i+z_i$. Заметим, что $v_i\geq w_i+x_i$. Тогда
\be
\label{R A_i - A A_i-1}
R(A_i)-R(A_{i-1}) \leq u_i-w_i-x_i\leq n(A_0)-3w_i+z_i-x_i.
\ee
На псевдопрямой $l_i$ находится $z_i+x_i$ точек пересечения, не принадлежащих псевдопрямым из $A_0$.  Поэтому
\be
\label{f a_i - f a_i-1}
f(A_i)-f(A_{i-1})=n(A_0)-w_i+z_i+x_i.
\ee
Сложив формулы (\ref{R A_i - A A_i-1}) и (\ref{f a_i - f a_i-1}) по всем $i=1,2,\ldots,I$  получим
\be
\label{R A i - R A 0}
R(A_i)-R(A_0)\leq In(A_0)-3\sum_{i=1}^{I}w_i+\sum_{i=1}^{I}z_i-\sum_{i=1}^{I}x_i. \ee

\be
\label{f A i - f A 0}
f(A_i)-f(A_0)= In(A_0)-\sum_{i=1}^{I}w_i+\sum_{i=1}^{I}z_i+\sum_{i=1}^{I}x_i.
\ee
Выразим $\sum_{i=1}^{I}w_i$ из неравенства (\ref{R A i - R A 0}) и подставим в равенство (\ref{f A i - f A 0}), учитывая $R(A_i)\geq 0$:
\be
\label{f A i ine}
f(A_i)\geq f(A_0)+\frac 23In(A_0)-\frac{R(A_0)}3+\frac 43\sum_{i=1}^{I}x_i+\frac 23\sum_{i=1}^{I}z_i.
\ee
Теперь рассмотрим два случая по отдельности.

{\it Случай (а).} Псевдопрямая $PQ$ не принадлежит семейству $A$. Тогда, подставляя в неравенство (\ref{f A i ine}) параметры
\bez
n(A_0)=2k, \ I=n-2k, \ f(A_0)=k^2+2k-1, \ R(A_0)=(k-1)^2+2, \ x_i\geq 0, \ z_i\geq 0,
\eez
получаем неравенство
\bez
f(A_i)\geq \frac{4kn-6k^2+8k-6}3.
\eez
Осталось заметить, что $f(A)=f(A_i)$ и что неравенство $\frac{4kn-6k^2+8k-6}3\geq (k+1)(n-k)$ равносильно неравенству $(k-3)(n-(3k-2))\geq 0$, которое выполняется при данных условиях на $n$ и $k$. Поэтому в случае (а) получаем $f\geq (k+1)(n-k)$.

{\it Случай (б).} Псевдопрямая $PQ$ принадлежит семейству $A$. Обозначим через $b_j$ количество точек пересечения псевдопрямых семейства $A$ кратности $j$, лежащих на псевдопрямой $PQ$ и отличных от точек $P$ и $Q$. Тогда количество псевдопрямых семейства $A$, не проходящих через точки $P$ и $Q$, равно $I=\sum_{j=2}^{k}(j-1)b_j=n-2k+1$. Из определения чисел $x_i$ и $z_i$ следует
\be
\label{z b x ine}
\sum_{i=1}^{I}z_i\geq b_3+\ldots+b_k \quad \text{и} \quad \sum_{i=1}^{I}x_i\geq \sum_{j=3}^{k}(j-3)b_j.
\ee
 Вычислим параметры семейства $A_0:$
\be
\label{A_0 n f }
n(A_0)=2k-2,  \ f(A_0)=k^2-2, \ R(A_0)=(k-2)^2+2.
\ee
Подставим (\ref{z b x ine}), (\ref{A_0 n f }) и формулу $f(A)-f(t_i)=2+\sum_{j=2}^{k}b_j$ в неравенство (\ref{f A i ine}):
\be
\label{f A b ine}
f(A)\geq \frac{(4k-4)n-6k^2+16k-10}3+b_2+\sum_{j=3}^{k}b_j\left(\frac{4j-7}3\right).
\ee
Так как
\bez
\sum_{j=2}^{k}(j-1)b_j=n-2k+1 \quad \text{ и} \quad \frac{4j-7}3\geq \frac 56(j-1)  \quad \text{ при} \quad j\geq 3,
\eez
то
\bez
b_2+\sum_{j=3}^{k}b_j\left(\frac{4j-7}3\right)\geq \frac 56(n-2k+1).
\eez
Учитывая последнее неравенство, преобразуем (\ref{f A b ine}) к виду
\bez
f(A)\geq \frac{(8k-3)n-12k^2+22k-15}6.
\eez
Осталось доказать неравенство
\bez
\frac{(8k-3)n-12k^2+22k-15}6\geq (k+1)(n-k) \quad \Leftrightarrow \quad  n(2k-9)\geq 6k^2-28k+15,
\eez
которое имеет место при $n\geq \frac{k^2+k}2+3\geq 3k+3$ и при $k\geq 6.$ Случай (б) разобран.
\ep

\end{document}